\documentclass{article}
\usepackage{amssymb,amsfonts,amsmath,amsthm,enumerate}
\usepackage[numeric]{amsrefs}

\def \a{\alpha}
\def \b{\beta}
\def \C{{\mathbb C}}
\def \ds{\displaystyle}

\def \f{\frac}
\def \t{\theta}
\def \z{\zeta}

\def \vp{\varphi}

\def \G{\Gamma}

\def \sgn{\operatorname{sgn}}
\def \tr{\operatorname{tr}}

\def \R{{\mathbb R}}
\def \Re{\operatorname{Re}}

\def \Z{{\mathbb Z}}
\def \l{\left}
\def \r{\right}
\def \hat{\widehat}
\def \tilde{\widetilde}
\def \z{\zeta}

\newtheorem{theorem}{Theorem}[section]

\newtheorem{corollary}[theorem]{Corollary}
\newtheorem{proposition}[theorem]{Proposition}
\numberwithin{equation}{section} \theoremstyle{definition}

\newtheorem*{remark}{Remark}

\begin{document}

\title{Triple mean values of Witten $L$-functions}
\author{Shin-ya Koyama\footnote{Department of Biomedical Engineering, Toyo University,
2100 Kujirai, Kawagoe, Saitama, 350-8585, Japan.} \ \& Nobushige Kurokawa\footnote{Department of Mathematics, Tokyo Institute of Technology, 
Oh-okayama, Meguro-ku, Tokyo, 152-8551, Japan.}}

\maketitle

\begin{abstract}
Mean values of Witten $L$-functions in the ``character'' aspect are
investigated. After giving a general formula for mean values with
the first and the second power, 
we explicitly calculate the cubic moment for $SU(2)$.
\end{abstract}

\section{Introduction}
Study of mean values of zeta functions is one of central subjects
in number theory.
The most frequently investigated problem would be the case for
the mean value of the $2k$-th power of zeta functions
in the $t$-aspect along the critical line  (the $2k$-th moment).
For instance,
much work has been done towards the conjectural asymptotic
for the Riemann zeta function $\zeta(s)$:
\begin{equation}\label{2k}
\f1{2T}\int_{-T}^T\z\l(\f12+it\r)^{2k}dt
\sim c_k(\log T)^{k^2}\qquad (T\to\infty).
\end{equation}
Analogous problems exist for various zeta functions for
more general $s\in\C$ in more general aspects.

We find a tendency that
the higher $k$ is, the more difficult the problem is.
For example, the value $c_k$ in \eqref{2k} 
as well as its proof is known only for $k=1,2$,
which are classical results by Hardy-Littlewood \cite{HL}
and Ingham \cite{Ingham}. The conjectural values of $c_k$ are known
only for $k=3,4$ (\cite{CGhosh} \cite{CGonek}).
A general form of $c_k$ is proposed by Keating-Snaith \cite{KS}
and Brezin-Hikami \cite{BH} under assuming the analogy
between $\zeta(s)$ and the characteristic polynomial of random matrices.

We also observe that any odd power moment is very hard to treat.
As far as the authors know, 
the only  successful case where the cubic moment was studied is
the work by Conrey and Iwaniec \cite{CI}.

The goal of this paper is to deal with the mean values of the third power of
Witten $L$-functions in the ``character'' aspect in $\Re(s)>0$.
For a compact semisimple Lie group $G$, Witten (\cite{Witten}) 
defined a zeta function
from the partition function of a quantum system as follows:
\begin{equation}\label{witten}
\zeta_G(s)=\sum_{\rho\in\hat G}(\dim\rho)^{-s},
\end{equation}
where $\hat G$ denotes the set of equivalence classes of
irreducible unitary representations of $G$.
It is known that \eqref{witten} is absolutely convergent 
if $\Re(s)$ is sufficiently large (\cite{Bd},\cite{LL}), 
and that $\zeta_G(s)$ is meromorphic in $s\in\C$.

In case of $G=SU(2)$, it holds that
$$
\hat G=\{\mathrm{Sym}^{n-1}\,|\,n=1,2,3,...\},
$$
where $\mathrm{Sym}^{n-1}$ is the $n$-dimensional
symmetric power representation defined by
$$
\mathrm{Sym}^{n-1}:\ SU(2)\ni g\longmapsto
\begin{pmatrix}
\a^{n-1} & & & \\
 &\a^{n-2}\b& & \\
&&\ddots&\\
&&&\b^{n-1}
\end{pmatrix}\in GL(n,\C)
$$
with $\a$ and $\b$ being the eigenvalues of $g$. Hence
$$
\zeta_{SU(2)}(s)=\sum_{n=1}^\infty n^{-s}=\zeta(s).
$$
In this sense the Witten zeta function is a generalization (a deformation)
of the Riemann zeta function.

We also define the Witten $L$-functions after Kurokawa-Ochiai \cite{KO}
by attaching ``characters'' which suitably twist $\rho\in\hat G$.
Any fixed element $g\in G$ plays the role.
Namely, we regard the following map as a ``character'' of $\hat G$:
$$
g:\ \hat G\ni\rho\longmapsto\f{\tr(\rho(g))}{\deg\rho}\in\C.
$$
Here $\f1{\deg\rho}$ is a normalization factor so that
$g(\rho)=1$ for the unit element $g\in G$.
In this manner the Witten $L$-function of $G$ with a twist given by
$g\in G$ is defined by
\begin{equation}\label{wittenL}
\zeta_G(s,g)=\sum_{\rho\in\hat G}\f{\tr(\rho(g))}{\deg\rho}(\deg\rho)^{-s}.
\end{equation}
This is also absolutely convergent for $s\in\C$ with $\Re(s)$ 
sufficiently large.
It is easy to see that $\zeta_G(s,g)$ depends only on the 
conjugacy class of $g$ in $G$.

The chief concern in this paper is to study the mean values
\begin{equation}\label{Z}
Z_G^m(s):=\int_G\z_G(s,g)^mdg
\end{equation}
with $dg$ the normalized Haar measure of $G$. 
Such a problem was first studied in \cite{KO} \S2.8.

In the next section we study \eqref{Z} for $m=1,2$ and $\Re(s)>1$.
In the final section we specialize the group as $G=SU(2)$,
and calculate \eqref{Z} for $m=3$.
We also show that the function $Z_{SU(2)}^3(s)$ is meromorphic in $s\in\C$.

\section{Preliminary results}
\begin{proposition}[The mean value: the first power moment]\label{2.1}
For any compact semisimple Lie group $G$, it holds that
$$
Z_G^1(s):=\int_G \z_G(s,g)dt=1,
$$
where $\Re(s)$ is sufficiently large.
\end{proposition}

\proof
Let $\Re(s)$ be large enough so that
the series \eqref{wittenL} is absolutely convergent.
Then
\begin{align}
\int_G \z_G(s,g)dt
&=\int_G \sum_{\rho\in\hat G}\f{\tr(\rho(g))}{\deg\rho}(\deg\rho)^{-s}dg\nonumber\\
&=\sum_{\rho\in\hat G}\f1{(\deg\rho)^{s+1}}\int_G \tr(\rho(g))dg.
\label{sum}
\end{align}
Here we appeal to the orthogonal relation of characters.
For any $\rho,$ $\rho'\in\hat G$, it holds that
\begin{equation}\label{orth}
\int_G \tr(\rho(g))\overline{\tr(\rho'(g))}dt
=\begin{cases}
1 &(\rho=\rho')\\
0 &(\text{otherwise}).
\end{cases}
\end{equation}
Putting $\rho'$ to be the identity representation,
we find that
$$
\int_G \tr(\rho(g))dg
=\begin{cases}
1 &(\rho=1)\\
0 &(\text{otherwise}).
\end{cases}
$$
Hence all terms in the sum in \eqref{sum} are zero except for $\rho=1$.
The proposition follows.
\qed

\begin{remark}
\begin{enumerate}[(1)]
\item
The assumption on $\Re(s)$ is necessary. Indeed from the known fact that
$\zeta_{SU(2)}(-2,g)=0$ for all $g\in SU(2)$
(\cite{KO} Theorem1), we obviously see that $\int_{SU(2)}\z(-2,g)dg=0$.
We also refer to \cite{GJK} on vanishing of Witten zeta functions at $s=-2$
for the compact $p$-adic Lie groups.
(See also \cite{MinZeros}, \cite{MinMeanvalue} for some generalizations.)

\item
The preceding remark shows that a meromorphic continuation of
\eqref{Z} is different from the integral of the meromorphically continued
$\z_G(s,g)^m$. Thus analyticity of $Z_G^m(s)$ is nontrivial, in the sense
that it is not an immediate consequence from that of $\z_G(s,g)$.
\end{enumerate}
\end{remark}

In the next proposition
we assume that the group $G$ satisfies the following condition:
\begin{equation}\label{*}\tag{$\ast$}
\tr(\rho(g))\in\R\qquad(\forall\rho\in\hat G,\ \forall g\in G).
\end{equation}
The special unitary group $G=SU(2)$ is an example of such $G$,
as shown in Proposition \ref{3.1} below.

\begin{proposition}[The double mean value: the square moment]\label{2.2}
For any compact semisimple Lie group $G$ satisfying the condition \eqref{*},
it holds that
$$
Z_G^2(s):=\int_G \z_G(s,g)^2dt=\zeta_G(2s+2),
$$
where $\Re(s)$ is sufficiently large.
\end{proposition}

\proof
Let $\Re(s)$ be large enough so that
the series \eqref{wittenL} is absolutely convergent.
Then
\begin{align*}
\int_G \z_G(s,g)^2dt
&=\int_G \sum_{\rho_1,\rho_2\in\hat G}
\f{\tr(\rho_1(g))\tr(\rho_2(g))}{(\deg\rho_1)^{s+1}(\deg\rho_2)^{s+1}}dg\\
&=\sum_{\rho_1,\rho_2\in\hat G}\f1{(\deg\rho_1)^{s+1}(\deg\rho_2)^{s+1}}
\int_G\tr(\rho_1(g))\tr(\rho_2(g))dg.
\end{align*}
By the orthogonal relation \eqref{orth},
all terms in the sum are zero except for $\rho_1=\rho_2$.
Therefore if we put $\rho_1=\rho_2=\rho$, it holds that
$$
\int_G \z_G(s,g)^2dt
=\sum_{\rho\in\hat G}\f1{(\deg\rho)^{2s+2}}=\z_G(2s+2).
$$
\qed

When $G$ is a finite group, we have
$$
Z_G^m(s)=\f1{|G|}\sum_{g\in G}\z_G(s,g)^m.
$$

The following proposition is an example of calculation of
mean values for general $m\ge1$.

\begin{proposition}[Mean values $Z_{S_3}^m(s)$]\label{sym}
Let $G=S_3$ be the symmetric group of degree 3.
Then it holds that
$$
Z_{S_3}^m(s)=\f{(2+2^{-s})^m+2(2-2^{-s-1})^m}6.
$$
Especially,
\begin{align*}
Z_{S_3}^1(s)&=1,\\
Z_{S_3}^2(s)&=2+2^{-2s-2}=\z_{S_3}(2s+2).
\end{align*}
\end{proposition}

\proof
The Witten $L$-function $\z_G(s,g)$ depends only on the conjugacy class of $g$.
Conjugacy classes of symmetric groups are classified by the cycle type.
The elements in $S_3$ consist of two cyclic permutations of order three,
three transpositions of order two, and the identity. Hence
\begin{align*}
Z_{S_3}^m(s)
&=\f1{6}\l(2\z_{S_3}(s,(1\,2\,3))^m+3\z_{S_3}(s,(1\,2))^m+\z_{S_3}(s,(1))^m\r).
\end{align*}
Now we will calculate $\z_{S_3}(s,g)$ for each $g\in S_3$.
We have $\hat{S_3}=\{\rho_1,\rho_2,\rho_3\}$ with
$\rho_1$ the trivial representation, $\rho_2$ the signature, 
and $\rho_3$ the unique two dimensional irreducible representation defined by
$$
\rho_3((1\,2))=\begin{pmatrix}0&1\\1&0\end{pmatrix},\qquad
\rho_3((1\,2\,3))=\begin{pmatrix}e^{\f23\pi i}&0\\0&e^{-\f23\pi i}\end{pmatrix},
$$
whose traces are 0 and $-1$, respectively.
Then each $L$-function is given by
\begin{align*}
\z_{S_3}(s,(1))&=1+1+\f2{2^{1+s}}=2+2^{-s},\\
\z_{S_3}(s,(1\,2))&=1-1+0\cdot2^{-s}=0,\\
\z_{S_3}(s,(1\,2\,3))&=1+1+\f{-1}{2^{1+s}}=2-2^{-s-1}.
\end{align*}
Therefore we get the conclusion.
\qed

The characters attached to
Witten $L$-functions are generalized to 
convolutions of $m$ times of twists $(m=1,2,3,...)$:
$$
\z_G(s,(g_1,...,g_m))=
\sum_{\rho\in\hat G}\f{\tr(\rho(g_1))\cdots\tr(\rho(g_m))}{(\deg\rho)^{s+m}}.
$$
We can consider another type of mean values as
\begin{equation}\label{another}
\tilde Z_G^m(s):=\int_G\z_G(s,(\underbrace{g,...,g}_m))dg
=\int_G\sum_{\rho\in\hat G}\f{(\tr(\rho(g)))^m}{(\deg\rho)^{s+m}}dg.
\end{equation}
This is calculated for $G=S_3$ in the following proposition.

\begin{proposition}
Let $S_3$ be the symmetric group of degree 3.
Then it holds that
$$
\tilde Z_{S_3}^m(s)=\int_{S_3}\z_{S_3}(s,(\underbrace{g,...,g}_m))dg
=\f{3+(-1)^m}2+\f{2^{-1}+(-1)^m}3\cdot 2^{-s-m}.
$$
\end{proposition}

\proof
We again use 
the classification of the conjugacy classes of $S_3$ and the explicit form
of elements in $\hat{S_3}$ given in the proof of the preceding proposition. 
It holds that
\begin{align*}
\lefteqn{\int_{S_3}\z_{S_3}(s,(\underbrace{g,...,g}_m))dg}\\
&=\f16\sum_{\rho\in\hat{S_3}}\l(
\f{2\tr(\rho(1\,2\,3))}{(\deg\rho)^{s+m}}+
\f{3\tr(\rho(1\,2))}{(\deg\rho)^{s+m}}+
\f{1}{(\deg\rho)^{s+m}}\r)\\
&=\f16\l(\underbrace{2+3+1}_{\rho_1}+
\underbrace{2+3(-1)^m+1}_{\rho_2}+
\underbrace{\f{2(-1)^m}{2^{s+m}}+\f{3\cdot 0}{2^{s+m}}+\f1{2^{s+m}}}_{\rho_3}\r)\\
&=\f{3+(-1)^m}2+\f{2^{-1}+(-1)^m}32^{-s-m}.
\end{align*}
\qed

\section{Triple mean values: the cubic moment}
\begin{proposition}\label{3.1}
The group $G=SU(2)$ satisfies the condition \eqref{*}.
\end{proposition}

\proof
We put the eigenvalues of $g\in SU(2)$ as $e^{\pm i\t}$.
The conjugacy classes of $G=SU(2)$ are parametrized by $\t\in[0,\pi]$.
Since the Witten $L$-function depends only on the conjugacy class of $g\in G$,
we denote $\zeta_G(s,g)=\zeta_G(s,[\t])$.
Then we compute
\begin{align}
\tr(\mathrm{Sym}^{n-1}(g))
&=\tr\l(\mathrm{Sym}^{n-1}\l(\begin{pmatrix}
e^{i\t}&\\&e^{-i\t}\end{pmatrix}\r)\r)\nonumber\\
&=\tr\begin{pmatrix}
e^{i(n-1)\t} & & & \\
 &e^{i(n-3)\t}& & \\
&&\ddots&\\
&&&e^{-i(n-1)\t}
\end{pmatrix}\nonumber\\
&=\begin{cases}
n & (\t=0)\\
\dfrac{\sin(n\t)}{\sin\t} & (0<\t<\pi)\\
(-1)^{n-1}n & (\t=\pi).\label{trsym}
\end{cases}
\end{align}
\qed

From \eqref{trsym}
we obtain the following corollary immediately.

\begin{corollary}[The explicit form of $\z_{SU(2)}(s,\lbrack\t\rbrack)$]
We have
$$
\z_{SU(2)}(s,[\t])=\begin{cases}
\zeta(s) & (\t=0)\\
\ds\sum_{n=1}^\infty\f{\sin(n\t)}{\sin\t}n^{-s-1} & (0<\t<\pi)\\
\ds\sum_{n=1}^\infty\f{(-1)^{n-1}}{n^s}=(1-2^{1-s})\z(s) & (\t=\pi).
\end{cases}
$$
\end{corollary}

Before calculationg the triple mean value, we are introducing
preliminary calculations which are good for general $m$.
For $G=SU(2)$ and $\Re(s)>1$, we compute by putting 
$\rho_j=\mathrm{Sym}^{n_j-1}$ that
\begin{align*}
\lefteqn{\int_G \z_G(s,g)^mdg}\\
&=\sum_{\rho_1,...,\rho_m\in\hat G}
\f1{((\deg\rho_1)\cdots(\deg\rho_m))^{s+1}}
\int_G\tr(\rho_1(g))\cdots\tr(\rho_m(g))dg\\
&=\sum_{n_1,...,n_m\ge1}
\f1{(n_1\cdots n_m)^{s+1}}
\int_G\tr(\mathrm{Sym}^{n_1-1}(g))\cdots\tr(\mathrm{Sym}^{n_m-1}(g))dg.
\end{align*}
The last integrand depends only on the conjugacy class of $g$.
We use the notation $[\t]$ defined in the proof of Proposition \ref{3.1}.
Then the last integral equals by \eqref{trsym}
\begin{align}
\lefteqn{
\int_0^\pi\tr\l(\mathrm{Sym}^{n_1-1}([\t])\r)
\cdots\tr\l(\mathrm{Sym}^{n_m-1}([\t])\r)
\f2\pi \sin^2\t dt}\nonumber\\
&=\int_0^\pi\f{\sin(n_1\t)}{\sin\t}\cdots\f{\sin(n_m\t)}{\sin\t}
\f2\pi \sin^2\t dt\nonumber\\
&=\f2\pi \int_0^\pi\f{\sin(n_1\t)\cdots\sin(n_m\t)}{\sin^{m-2}\t}d\t.
\label{cn1nm}
\end{align}
Putting it by $c(n_1,...,n_m)$, we have for $\Re(s)>1$ that
\begin{align}\label{c1...m}
Z_G^m(s):=\int_G \z_G(s,g)^mdg
&=\sum_{n_1,...,n_m\ge1}\f{c(n_1,...,n_m)}{(n_1\cdots n_m)^{s+1}}.
\end{align}

\begin{theorem}[Triple mean values: the cubic moment]\label{cubic}
It holds for $G=SU(2)$ in $\Re(s)>0$ that
\begin{align}
Z_G^3(s):=&\int_G \z_G(s,g)^3dg\nonumber\\
=&\sum_{m_1,m_2,m_3\ge0}\l((m_1+m_2+1)(m_2+m_3+1)(m_3+m_1+1)\r)^{-s-1}.
\label{dirichlet}
\end{align}
This function is meromorphically continued to all $s\in\C$.
\end{theorem}

\proof
We compute \eqref{c1...m} for $m=3$.
By transforming the product of the sine function into sums, we have
\begin{align}
\lefteqn{c(n_1,n_2,n_3)}\nonumber\\
&=\f2\pi\int_0^\pi\f{\sin(n_1\t)\sin(n_2\t)\sin(n_3\t)}{\sin\t}d\t\nonumber\\
&=\f1{2\pi}\int_0^\pi\l(
\f{\sin((n_1+n_2-n_3)\t)}{\sin\t}+
\f{\sin((n_2+n_3-n_1)\t)}{\sin\t}\r.\nonumber\\
&\qquad\qquad\qquad \l.+
\f{\sin((n_3+n_1-n_2)\t)}{\sin\t}-
\f{\sin((n_1+n_2+n_3)\t)}{\sin\t}\r)d\t\nonumber\\
&=\f{A(n_1+n_2-n_3)+A(n_2+n_3-n_1)+A(n_3+n_1-n_2)-A(n_1+n_2+n_3)}2
\label{c123}
\end{align}
with
$$
A(n):=\f1\pi\int_0^\pi\f{\sin(n\t)}{\sin\t}d\t.
$$

We first compute $A(n)$ for $n\ge1$. It holds that
\begin{align*}
A(n)
&=\f1\pi\int_0^\pi\f{e^{in\t}-e^{-in\t}}{e^{i\t}-e^{-i\t}}d\t\\
&=\f1\pi\int_0^\pi\l(e^{i(n-1)\t}+e^{i(n-3)\t}+\cdots+e^{-i(n-1)\t}\r)d\t\\
&=\begin{cases}
1 &(n:\text{ odd})\\
0 &(n:\text{ even}),
\end{cases}
\end{align*}
because the integrand contains the constant term ``1''
if and only if $n$ is odd, which contribute 1 to $A(n)$,
and the integral vanishes for all other terms.
Since $A(n)$ is an odd function in $n$, we eventually have for $n\in\Z$ that
$$
A(n)
=\begin{cases}
\sgn(n)&(n:\text{ odd})\\
0 &(n:\text{ even}).
\end{cases}
$$

Next we calculate \eqref{c123}.
When $n_1+n_2+n_3$ is even, all of
$n_1+n_2-n_3$, $n_2+n_3-n_1$ and $n_3+n_1-n_2$ are even,
and thus $c(n_1,n_2,n_3)=0$.
Assume that $n_1+n_2+n_3$ is odd. Then all of
$n_1+n_2-n_3$, $n_2+n_3-n_1$ and $n_3+n_1-n_2$ are odd, and
$$
c(n_1,n_2,n_3)
=\f{\sgn(n_1+n_2-n_3)+\sgn(n_2+n_3-n_1)+\sgn(n_3+n_1-n_2)-1}2.
$$
Here we put the following condition as ($\ast\ast$):
\begin{multline}\label{**}\tag{$\ast\ast$}
\text{All of $n_1+n_2-n_3$, $n_2+n_3-n_1$ and $n_3+n_1-n_2$ are positive,}\\
\text{and $n_1+n_2+n_3$ is odd.}
\end{multline}
When the triple $(n_1,n_2,n_3)$ satisfies ($\ast\ast$),
we have $c(n_1,n_2,n_3)=1$.
Assume that the triple $(n_1,n_2,n_3)$ does not satisfy ($\ast\ast$).
Then it is easy to see that only one of 
$n_1+n_2-n_3$, $n_2+n_3-n_1$ and $n_3+n_1-n_2$ is negative or zero.
But it cannot be zero, because it is odd by assumption.
So one of $n_1+n_2-n_3$, $n_2+n_3-n_1$ and $n_3+n_1-n_2$ is negative,
and the other two are positive.
Hence we conclude that $c(n_1,n_2,n_3)=\f{1+1-1-1}2=0$, when
the triple $(n_1,n_2,n_3)$ does not satisfy ($\ast\ast$).
Therefore we successfully have the final form of the coefficients as
\begin{equation}\label{10}
c(n_1,n_2,n_3)
=\begin{cases}1&(\ast\ast)\\0&(\text{otherwise}).\end{cases}
\end{equation}
By \eqref{c1...m} the triple mean value is
\begin{equation}\label{n123}
Z_G^3(s)=\int_G \z_G(s,g)^3dg
=\sum_{\genfrac{}{}{0pt}{1}{n_1,n_2,n_3\ge1}{(\ast\ast)}}\f1{(n_1n_2n_3)^{s+1}}.
\end{equation}
We have proved \eqref{n123} for $\Re(s)>1$, but now we see that
the right hand side of \eqref{n123} is absolutely convergent in 
$\Re(s)>0$, because its absolute value is $O\l(\zeta(s)^3\r)$. 
So $Z_G^3(s)$ is analytically continued to $\Re(s)>0$
by \eqref{n123}.

In what follows we rewrite \eqref{n123} to a simpler form.
Put
$$
m_3:=\f{n_1+n_2-n_3-1}2,\quad
m_1:=\f{n_2+n_3-n_1-1}2,\quad
m_2:=\f{n_3+n_1-n_2-1}2.
$$
Then there is one-to-one correspondence between
the set of all triples $(n_1,n_2,n_3)$ with ($\ast\ast$)
and the set of all triples $(m_1,m_2,m_3)\in(\Z_{\ge0})^3$.
The inverse correspondence
$$
n_1=m_2+m_3-1,\quad
n_2=m_3+m_1-1,\quad
n_3=m_1+m_2-1
$$
leads to the conclusion.

The meromorphic continuation was generally shown by Mellin \cite{Mellin}.
\qed

By this theorem the first few terms of $Z_G^3(s)$ turns to be as follows:
$$
Z_G^3(s)=
1+\f3{4^s}+\f3{9^s}+\f3{12^s}+\f3{16^s}+\f6{24^s}+\f3{25^s}+\f1{27^s}+\cdots.
$$

We can also compute the other type \eqref{another} of triple mean value
for $G=SU(2)$ and $m=3$.
\begin{theorem}
It holds for $\Re(s)>-2$ that
$$
\tilde Z_{SU(2)}^3(s)=\int_{SU(2)}\z_{SU(2)}(s,(g,g,g))dg=(1-2^{-s-3})\z(s+3).
$$
In particular, the function $\tilde Z_{SU(2)}^3(s)$ is meromorphic in $\C$.
\end{theorem}

\proof
We compute
\begin{align}
\int_{SU(2)}\z_{SU(2)}(s,(g,g,g))dg
&=\sum_{n=1}^\infty\f1{n^{s+3}}\int_0^\pi\l(\f{\sin n\t}{\sin\t}\r)^3
\f2\pi\sin^2\t d\t\nonumber\\
&=\sum_{n=1}^\infty\f{c(n,n,n)}{n^{s+3}}\label{3.7},
\end{align}
where $c(n_1,n_2,n_3)$ is defined by \eqref{cn1nm}.
By applying \eqref{10} with $m=3$ and $n_1=n_2=n_3$, we find that
$$
c(n,n,n)
=\begin{cases}1&(n:\ \text{odd})\\0&(\text{otherwise}).\end{cases}
$$
Then \eqref{3.7} is absolutely convergent in $\Re(s)>-2$ and it holds that
$$
\int_{SU(2)}\z_{SU(2)}(s,(g,g,g))dg
=\sum_{\genfrac{}{}{0pt}{1}{n\ge1}{\text{odd}}}\f1{n^{s+3}}
=(1-2^{-s-3})\z(s+3).
$$
\qed

\begin{remark}[Quadruple case]
The case of $m=4$ is also calculated as follows.
Starting from the identity
\begin{align*}
\tilde Z_{SU(2)}^4(s)
&=\int_{SU(2)}\z_{SU(2)}(s,(g,g,g,g))dg\\
&=\sum_{n=1}^\infty\f1{n^{s+4}}\int_0^\pi\l(\f{\sin n\t}{\sin\t}\r)^4
\f2\pi\sin^2\t d\t\\
&=\sum_{n=1}^\infty\f1{n^{s+4}}\cdot\f2\pi\int_0^\pi\f{\sin^4 n\t}{\sin^2\t}d\t,
\end{align*}
the last integral is calculated as follows
\begin{align*}
\lefteqn{\f2\pi\int_0^\pi\f{\sin^4 n\t}{\sin^2\t}d\t}\\
&=\f{-1}{2\pi}\int_0^\pi\l(
\f{e^{in\t}-e^{-in\t}}{e^{i\t}-e^{-i\t}}(e^{in\t}-e^{-in\t})\r)^2d\t\\
&=\f{-1}{2\pi}\int_0^\pi\l(
(e^{i(n-1)\t}+e^{i(n-3)\t}+\cdots+e^{-i(n-1)\t})(e^{in\t}-e^{-in\t})\r)^2d\t\\
&=\f{-1}{2\pi}\int_0^\pi\l((e^{i(2n-1)\t}+e^{i(2n-3)\t}+\cdots+e^{i\t})
-(e^{-i\t}+e^{-3i\t}+\cdots+e^{-i(2n-1)\t})\r)^2d\t\\
&=\f{-1}{2\pi}\int_0^\pi(-2n+[\text{nonconstant terms}])d\t=n.
\end{align*}
For the nonconstant terms are written as a linear combination of
$$e^{iN\t}+e^{-iN\t}=2\cos(N\t)\qquad (N\in\Z\setminus\{0\}),$$
whose integral vanishes as
$$
\int_0^\pi\cos(N\t)d\t=0\qquad (\forall N\in\Z\setminus\{0\}).
$$
Therefore we conclude that
\begin{equation}\label{quadruple}
\tilde Z_{SU(2)}^4(s)
=\sum_{n=1}^\infty\f1{n^{s+4}}\cdot n=\z(s+3).
\end{equation}
This is valid for $\Re(s)>-2$, and it shows that
$\tilde Z_{SU(2)}^4(s)$ has an analytic continuation to the entire plane
except for a simple pole at $s=-2$.

We can directly calculate that
$\int_{SU(2)}\z_{SU(2)}(-2,(g,g,g,g))dg=\infty$
by using the result of 
Min \cite{MinZeros} on the values $\z_{SU(2)}(-2,(g,g,g,g))$.
This happens to agree to our conclusion \eqref{quadruple} 
that $\tilde Z_{SU(2)}^4(s)$ has a pole at $s=-2$.
Although the problems are different as noted in Remark (2) after
Proposition 2.1, 
reasoning of this coincidence may be an interesting problem
for our future study.
\end{remark}

\section{Generalizations}
Let $G$ be as above, and $H$ be a subgroup of $G$.
Mean values are generalized to the average over $H$ as
$$
Z_{G,H}^m(s):=\int_H\z_G(s,h)^mdh
$$
and
$$
\tilde Z_{G,H}^m(s):=\int_H\z_G(s,(\underbrace{h,...,h}_m))dh.
$$
It is easy to see that
$Z_{G,\{1\}}^m(s)=\z_G(s)^m$ and
$Z_{G,G}^m(s)=Z_G^m(s)$.

\begin{theorem}
Let $G$ and $H$ be
a pair of compact semisimple Lie groups such
that $G=\underbrace{H\times\cdots\times H}_m$.
We regard $H$ as a subgroup of $G$ by diagonal embedding.
Then the following identities hold.
\begin{enumerate}[(1)]
\item
$
\z_G(s,(\underbrace{h,...,h}_m))=\z_H(s,h)^m.
$
\item
$
\tilde Z_{G,H}^m(s)=Z_{H}^m(s).
$
\end{enumerate}
\end{theorem}

\proof
(1)\ The map
$$
\underbrace{\hat H\times\cdots\times\hat H}_m
\ni(\rho_1,...,\rho_m)\longmapsto\rho_1\boxtimes\cdots\boxtimes\rho_m\in\hat G
$$
defined by
$$
(\rho_1\boxtimes\cdots\boxtimes\rho_m)(h_1,...,h_m)
:=\rho_1(h_1)\otimes\cdots\otimes\rho_m(h_m)
$$
is an isomorphism. It also holds that 
$$
\deg(\rho_1\boxtimes\cdots\boxtimes\rho_m)=(\deg\rho_1)\cdots(\deg\rho_m).
$$
Hence
\begin{align*}
\z_G(s,(\underbrace{h,...,h}_m))
&=\sum_{\rho\in\hat G}\f{\tr(\rho(h,...,h))}{(\deg\rho)^{s+1}}\\
&=\sum_{\rho_1,...,\rho_m\in\hat H}\f{\tr(\rho_1(h))\cdots\tr(\rho_m(h))}
{((\deg\rho_1)\cdots(\deg\rho_m))^{s+1}}\\
&=\l(\sum_{\rho\in\hat H}\f{\tr(\rho(h))}{(\deg\rho)^{s+1}}\r)^m
=\z_H(s,h)^m.
\end{align*}
(2)\ By (1), we have
$$
\tilde Z_{G,H}^m(s)=\int_H\z_G(s,(\underbrace{h,...,h}_m))dh
=\int_H\z_H(s,h)^mdh=Z_{H}^m(s).
$$
\qed

By Theorem \ref{cubic}, the following theorem is immediate.

\begin{theorem}\label{4.1}
Put $G=SU(2)\times SU(2)\times SU(2)$, and let $H=SU(2)$ be a
subgroup of $G$ embedded diagonally.
Then $\tilde Z_{G,H}^3(s)$ is explicitly expressed by the Dirichlet series
\eqref{dirichlet} in $\Re(s)>0$, and has a meromorphic
continuation to the whole plane $\C$.
\end{theorem}

\end{document}